# A New Theorem Relating the Tangent Secant Theorem to the Golden Ratio

M. N. Tarabishy

goodsamt@gmail.com

**Abstract**

The golden ratio is usually shrouded in mystique and mystery, however, showing its emergence from a familiar geometric setting makes it a more natural phenomenon. In this work, we present a new theorem connecting the Tangent Secant theorem to the golden ratio. Such a connection enhances the learning process of the theorem, and it demystifies the golden ratio in physics and other disciplines.

## 1. Introduction:

Science doesn't progress by new results only but also by finding new connections between known things. Our work here goes under the second category where we find a relation between two well known facts. The first is golden ratio which is a very interesting number that has been known since antiquity. The other known fact is the tangent secant theorem that has been known at least since Euclid, and while both of these facts are quite old, our search for a formal relationship between the two has turned up empty, and the closest we have found is in a comprehensive site about the golden ratio in geometry [1], and in this particular case, it came in relation to a right angled triangle of sides 3-4-5 and a circle used to get the golden ratio.

In the following sections we briefly introduce the golden ratio, then we introduce the tangent secant theorem, and after that the connection between the two is introduced in a form of a theorem.

## 2. The Golden Ratio:

The golden ratio φ is defined visually as a line shown in figure 1 that is divided into two segments a and b, such that the ratio of the large part to the small part is equal to the ratio of the sum of the two pieces to the large one, or:

$$\varphi = \frac{a}{b} = \frac{a+b}{a} = 1 + \frac{1}{\varphi} \qquad (1)$$

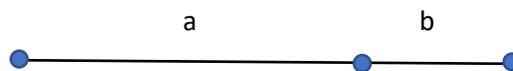

Figure 1. A line segment divided into two parts a, and b

Equation (1) can be written as a second order equation in φ:

$$\varphi^2 - \varphi - 1 = 0 \tag{2}$$

It is only appropriate here to use prof. Loh method [2] to solve this quadratic equation:

$\varphi_1 + \varphi_2 = (u +b) + (u-b) = 1 \rightarrow \varphi_{symmetry} = u = ½$

$\varphi_1 . \varphi_2 = (u+b) . (u-b) = -1 \rightarrow b^2 = u^2 +1 = 5/4 \rightarrow b = √5 /2 \rightarrow$ positive root:

$$\varphi = \frac{1 + \sqrt{5}}{2} = 1.6180 \ldots \tag{3}$$

Another way for representing the golden ratio that reflects its self-repeating nature is also obtained from equation (1) as:

$$\varphi = 1 + \frac{1}{\varphi} = 1 + \cfrac{1}{1 + \cfrac{1}{1 + \cfrac{1}{1 + \cfrac{1}{1 + \cfrac{1}{1 + \cdots}}}}} \tag{4}$$

### 2.1 The Golden Ratio and Fibonacci Sequence:

The famous sequence was popularized by the Italian mathematician Fibonacci's book: Liber Abaci in 1202. It is constructed by setting the next number as the sum of the two previous ones. Starting with the two numbers: 0, and 1, then, the next number is 0+1=1, and then the next is 1+1=2, then, the next one is 1+2=3, etc.…

$$0\ 1\ 1\ 2\ 3\ 5\ 8\ 13\ 21\ 34\ 55\ 89\ 144\ \ldots \tag{5}$$

Notice that:

$$\frac{21}{13} = 1.615, \quad \frac{34}{21} = 1.619, \quad \frac{89}{55} = 1.6181, \quad \ldots \rightarrow \varphi \tag{6}$$

So, the ratio of two adjoining numbers approaches the golden ratio φ as the numbers grow.

These numbers appear in so many natural settings like the number of petals in many flowers.

### 2.2 The Golden Ratio in Nature:

The golden ratio is a very interesting number on its own as it is the most irrational number, and this property is used by many plants for the efficient distribution of leaves and to fill a limited space with cells as can be seen in figure 2 where the green particles are arranged in golden ratio related spirals.

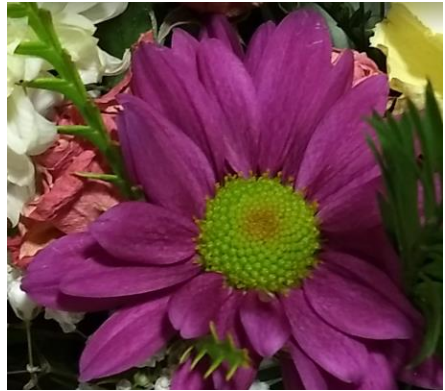

Figure 2. Arrangement of green cells in a flower

The website, Math is Fun [3], has a very interesting program that lets the user enter different rotation angles for arranging the cells and it shows how each choice affects the distribution of these cells, demonstrating the advantage of using the golden ratio.

The golden ratio has been known since ancient time, it can be observed in the great pyramid of Cheops that has a height, h = 146.515 m, and base, 2b = 230.363 m, [4].

The right triangle obtained by half the cross section of the pyramid is called the Egyptian (or Kepler triangle) that has an inclination of about 52° (51.82729°) and is shown in figure 3. The ratio s/b = 1.611 is the golden ratio within the measurement accuracy.

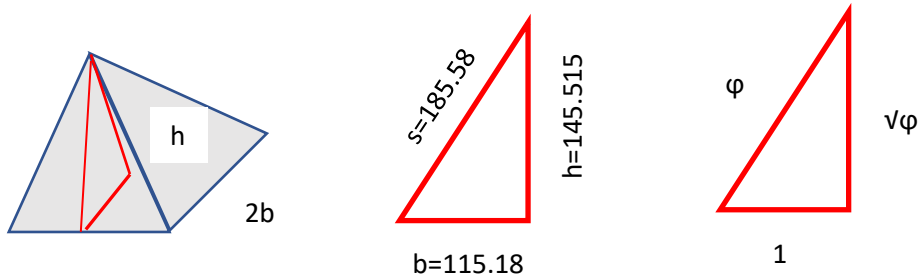

Figure 3. The dimensions of the great pyramid of Cheops and the Egyptian triangle

The Greeks wrote about the golden ratio and arguably used it in some of their structures, and some renaissance artists use it in their work. Recently, a study of the size of the human skulls [5] has shown golden ratio proportions. It has also appeared in symmetry signature of a nanoscale quantum system [6].

While it does appear in a number of interesting ways, its role in nature might have been sometimes exaggerated or at least not firmly proven like the claim that the ultimate beauty occurs when the body proportions are according to the golden ratio. However, even Professor Budd [7] who is critical of many claims about the golden ratio acknowledges its unique properties.

### 3. The Tangent-Secant Theorem:

The tangent-secant theorem is a very important theorem and is described by Euclid [8]. It relates the tangent and the secant emanating from the same point to a circle of radius r and centered at point w as shown in figure 4.

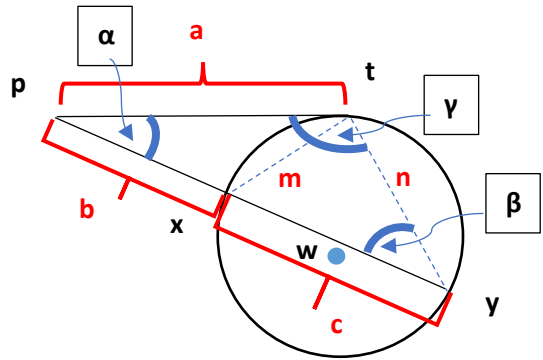

Figure 4. Tangent-Secant theorem diagram

If we denote the length of the tangent pt = a, and the outside secant px = b, the chord xy = c, and the secant s = py, then,

$$a^2 = b.(b+c) = b.s \qquad (7)$$

Equation (7) follows from the similarity of triangles: Δ pyt (the large triangle), and Δ pxt (the upper triangle) where the angle pyt = ptx = β , and angle pty = pxt = γ, and (b+c)/a = a/b = n/m.

### 4. New Theorem:

If we have a tangent and a secant to a circle from a point p as in figure 4, and if the ratio: tangent /chord = 1, then, the ratio of secant / tangent = φ, or,

$$If \ c = a, then, \qquad \frac{a}{b} = \frac{a+b}{a} = \varphi \qquad (8)$$

The opposite is true:

$$If \ \frac{a}{b} = \frac{c+b}{a} = \varphi, then, \qquad c = a \qquad (9)$$

Proof: From the tangent-secant theorem: a² = b . (b+c). If the chord c is equal to the tangent a, then we get: a² = b . (b+a), or, a/b = (a+b)/a = φ by definition.

For the opposite case, eqn. (9) gives (c/a) + (1/ φ) = φ, then, (c/a) = φ – ( φ-1) = 1.

Where: 1/φ = φ-1 according to eqn. (1)

### 4.1 Relationship Between Angles α and β:

To determine the conditions of the golden ratio occurrence, we examine the relationships that different triangles give. From the inscribed triangle, we get:

$$\frac{m}{\sin(\beta)} = \frac{n}{\sin(\alpha+\beta)} = \frac{c}{\sin(\gamma-\beta)} = 2r \qquad (10)$$

Where r is the circle radius, and we get:

$$\beta = \beta_1 = \frac{Arcsin\left(\frac{c}{2r}\right) - \alpha}{2} \qquad (11)$$

The upper triangle gives:

$$\frac{a}{\sin(\gamma)} = \frac{b}{\sin(\beta)} = \frac{m}{\sin(\alpha)} \qquad (12)$$

Resulting in:

$$\frac{a}{b} = \frac{\sin(\alpha+\beta)}{\sin(\beta)} = \varphi \qquad (13)$$

Also, we can get another expression:

$$\beta = \beta_2 = Arcsin\left[\sqrt{\frac{c}{2\,r\,\varphi}\sin(\alpha)}\,\right] \qquad (14)$$

Clearly, the relation between the angles α and β is nonlinear and is dependent on the chord to diameter ratio: c/2r (with values between 0 and 1). Equation (11), (13), and (14) can be used to find the corresponding angles.

Figures 5, 6, and 7 give examples of different angle α solutions for different c/2r ratio. Once we have the angle α, then we can get angles β, and γ.

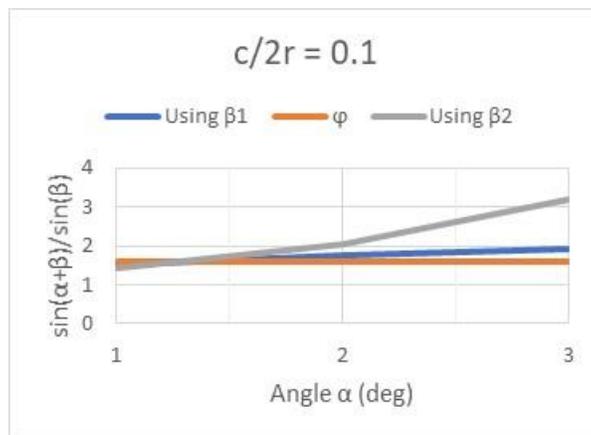

Figure 5. Angle α that corresponds to chord/diameter ratio (c/2r) = 0.1

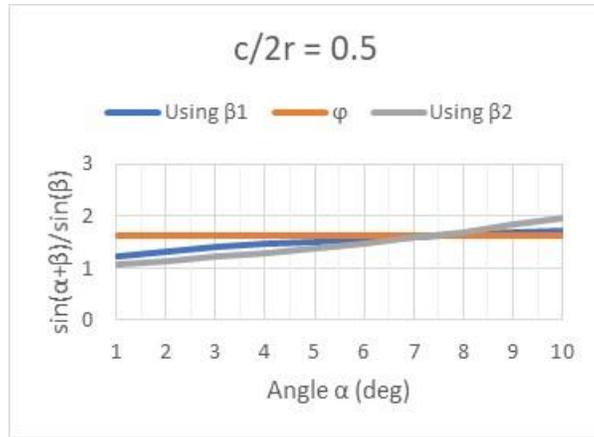

Figure 6.  Angle α that corresponds to chord/diameter ratio (c/2r) = 0.5

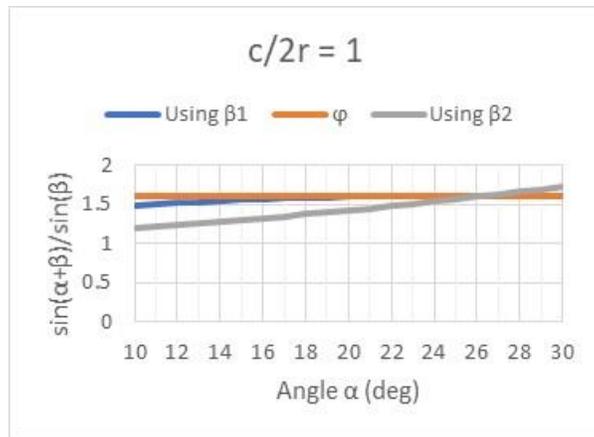

Figure 7.  Angle α that corresponds to chord/diameter ratio (c/2r) = 1

In the last case where the chord length is equal to 2r (the diameter) or c/2r = 1, as shown in figure 8, it is easy to calculate the angles geometrically (since the chord c is the diameter):

α = 26.565°, β = 31.717°, γ = 121.717°.

For the large triangle Δ pty, if the tangent a = T = 1, then, then, secant py= b+c = s = 1.618, n = 0.851.  For the upper triangle Δ ptx: a =1, px = 0.618, tx = m = 0.526.

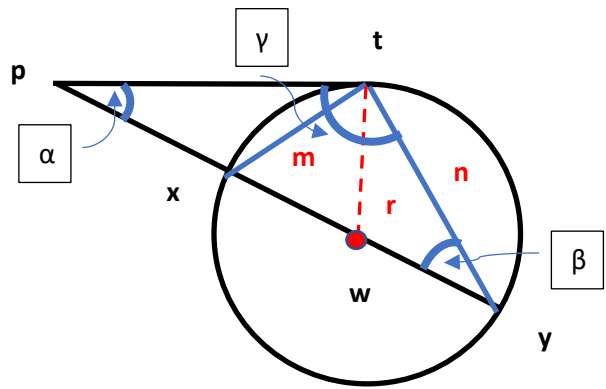

Figure 8. Example when the secant goes through the center of the circle

## 5. Discussion and Conclusion:

We have examined the golden ratio φ, its definition and its applications, then we reviewed the tangent secant theorem and then we showed the relationship between this theorem and the golden ratio by introducing a new theorem that states: If the tangent is equal to the chord, then, the ratio of the secant to the tangent is equal to ratio of the tangent to the outside secant and equals to the golden ratio φ. Then we showed examples of different values of c/2r: chord / diameter ratio when the secant/tangent arrangement is equal to the golden ratio. We found the angles graphically for three cases (c/2r = 0.1, 0.5, 1) using equations (11), (13), and (14). The new theorem helps in learning and remembering the tangent secant theorem and it also demystifies the golden ratio in physics and other disciplines.